\documentclass[11pt]{amsart}
\usepackage{amssymb,amsbsy,latexsym,amsmath,bbm,epsfig,psfrag,amsthm,mathrsfs}
\usepackage[swedish,english]{babel}
\usepackage{graphicx,tikz,esint}
\usepackage[T1]{fontenc}
\usepackage[latin1]{inputenc}
\usepackage{amsfonts}
\usepackage{amsmath}
\usepackage{amssymb}
\usepackage{epsf}
\usepackage{graphics}
\usepackage{graphicx}
\usepackage{latexsym}
\usepackage{psfrag}
\usepackage{relsize}
\usepackage{verbatim}
\usepackage{hyperref}
\usepackage{pgfplots}
\usepackage{float}
\usetikzlibrary{positioning}
\usetikzlibrary{arrows}
\usetikzlibrary{decorations.pathreplacing}
\usepackage[all]{xy}
\usepackage{url}

\theoremstyle{plain}

\newtheorem{Lem}{Lemma}
\numberwithin{Lem}{section}
\newtheorem{Prop}{Proposition}
\numberwithin{Prop}{section}
\newtheorem{Thm}{Theorem}
\numberwithin{Thm}{section}

\numberwithin{Cor}{section}

\numberwithin{Con}{section}

\theoremstyle{definition}
\newtheorem{Def}{Definition}
\numberwithin{Def}{section}

\numberwithin{hyp}{section}

\numberwithin{conj}{section}
\newtheorem{ex}{Example}
\numberwithin{ex}{section}

\newtheorem*{OP}{Open Problem}
 
\theoremstyle{remark}
\newtheorem{rem}{\bf{Remark}}
\numberwithin{rem}{section}

\numberwithin{equation}{section}

\DeclareMathOperator*{\tr}{tr\,}
\DeclareMathOperator*{\adj}{adj\,}
\DeclareMathOperator*{\Div}{div}
\DeclareMathOperator*{\Co}{co}
\DeclareMathOperator*{\Rco}{Rco}
\DeclareMathOperator*{\Int}{int}

\newcommand{\dv}{\partial}
\newcommand{\Om}{\Omega}

\newcommand{\eps}{\varepsilon}

\newcommand{\R}{{\mathbb R}}

\newcommand{\LL}{\mathcal{L}}

\begin{document}

\vspace{.4cm}

\title[]{  Generic ill-posedness of the energy-momentum equations \\ and differential inclusions}
\author[]{Erik Duse}
\date{}
\maketitle

\begin{abstract}
We show that the energy-momentum equations arising from inner variations whose Lagrangian satisfies a generic symmetry condition are ill-posed. This is done by proving that there exists a subclass of Lipschitz solutions that are also solutions to a differential inclusion into the orthogonal group and in particular these solutions can be nowhere $C^1$. We prove that these solutions are not stationary points if the Lagrangian $W$ is $C^1$ and strictly rank-one convex. In view of the Lipschitz regularity result of Iwaniec, Kovalev and Onninen for solution of the energy-momentum equation in dimension 2 we give a sufficient condition for the non-existence of a partial $C^1$ -regularity result even under the condition that the mappings satisfy a positive Jacobian determinant condition. Finally we consider a number of well-known functionals studied in nonlinear elasticity and geometric function theory and show that these do not satisfy this obstruction to partial regularity.
\end{abstract}

\section*{ Introduction}

Consider a functional 
\begin{align}\label{eq:F}
\mathcal{W}[u]=\int_{\Om}W(Du(x))dx
\end{align}

defined on mappings $u:\Om \subset \R^n\to \R^m$, where $\Om$ is an open set and $Du(x)$ denotes the total derivative of $u$ at $x$. Here, $\mathcal{W}[u]$ could for example be the energy of an elastic deformation of a hyperelastic material, in which case $n=m$, but there are other examples from physics in which $n\neq m$. In this paper we will however exclusively consider the case $n=m$. If $u\in C^1(\Om,\R^m)\cap C^0(\overline{\Om},\R^m)$ is a critical point of \eqref{eq:F}, then by considering \emph{outer variations} $u_\eps=u+\eps\phi$ for $\phi\in C^\infty_0(\Om,\R^m)$ one can show that $u$ solves the Euler-Lagrange equations in weak form 
\begin{align}\label{eq:Var}
\int_\Om \langle DW(Du(x)),D\phi(x)\rangle dx=0
\end{align}
for all $\phi\in C^\infty_0(\Om,\R^m)$, and if in addition $u\in C^2$, the $u$ also solves the Euler-Lagrange equations 
\begin{align}\label{eq:EL}
\text{div}\,DW(Du(x))=0
\end{align}
in the classical sense in $\Om$. Of course when considering weak solutions we can also consider much less regular functions $u$, so that $u$ is only Lipschitz continuous or $u$ is a Sobolev mapping for example. 

The outer variations are however not the only type of variations we may consider. Instead we could consider the \emph{inner-variations} instead. These are generated by a one-parameter family of diffeomorphisms $\phi_\eps\in C^\infty(\overline{\Om},\R^n)$ which in turn is generated by a smooth compactly supported vector field $\lambda$ so that 
\begin{align*}
\phi_\eps(x)=x+\eps\lambda(x). 
\end{align*}
The \emph{energy-momentum equations}, also called \emph{Noether's equations}, are then given by considering the variations $u_\eps=u(\phi_\eps)$ and setting
\begin{align*}
\frac{d}{d\eps}\bigg\vert_{\eps=0}\mathcal{W}(u_\eps)=0.
\end{align*}

The weak form of  \emph{energy-momentum equations} is given by
\begin{align}\label{eq:WeakEM}
\int_{\Om}\langle Du(x)^\ast DW(Du(x))-W(Du(x))I,D\lambda(x)\rangle dx=0
\end{align}
for all $\lambda \in C^\infty_0(\Om, \R^n)$. For a derivation see \cite[p. 147-150]{GH}.  By defining the \emph{energy-momentum tensor}
\begin{align}\label{eq:T}
T(x):=Du(x)^\ast DW(Du(x))-W(Du(x))I,
\end{align}
we see that $T(x)$ is divergence free in the sense of distributions. In the case of a $C^2$-solution there is a relation between the Euler-Lagrange equations and the energy-momentum equation given by the identity 
\begin{align}\label{eq:TEL}
\Div T(x)=Du(x)^\ast\Div DW(Du(x)).
\end{align}
Thus, if $u$ is a $C^2$-solution to the Euler-Lagrange equations it is also a solution to the energy-momentum equations. Conversely, if $u$ is a $C^2$-solution of the energy-momentum equations and $Du(x)^\ast$ is everywhere injective, then $u$ is also a solution to the Euler-Lagrange equations. 

In general however, and in particular in the case of vector valued mappings $u$, solutions of either \eqref{eq:EL} or \eqref{eq:TEL} are generically not $C^2$, and the weak form of the equations should be considered as independent conditions. In particular, there are weak solutions of the Euler-Lagrange equations associated to outer variations that are not weak solutions of \eqref{eq:TEL} and vice versa, see for example \cite{J}. Moreover, any \emph{strong local minimiser} of \eqref{eq:Var} satisfy both the weak Euler-Lagrange and the weak energy-momentum equations (provided $W$ satisfies some suitable structural conditions, see \cite[Thm. 2.4]{Ball}). For the definitions of weak and strong local minimiser we refer the reader to \cite[Ch. 4]{GH} and \cite{T}). 

Another conspicuous difference between the Euler-Lagrange and energy-momentum equations regards the formal determinedness of the equations. If $u:\Om\subset \R^n\to \R^m$ then the Euler-Lagrange equations are always formally determined whereas for the energy-momentum equations they are overdetermined if $n>m$, determined if $n=m$ and underdetermined if $n<m$.
Moreover, in nonlinear elasticity and geometric function theory for instance, we only want to minimise \eqref{eq:F} among mappings that are also homeomorphisms. In this case we want to have 
\begin{align}\label{eq:Jac}
J(x,u)=\det(Du(x))> 0\text{ a.e.}
\end{align}
in $\Om$. This point-wise constraint is a priori incompatible with the outer variations as they may violate the constraint. Therefore, it is unknown if the Euler-Lagrange equations hold or not even for minimisers. The inner variations however are compatible with the constraint \eqref{eq:Jac}, and for many natural material models in elasticity one can show that the weak energy momentum equations are satisfied (see \cite{Ball2}). 

These facts raise the natural question of what can be said about solutions of \eqref{eq:WeakEM}, for example under Dirichlet boundary conditions. Could these equations serve as a substitute for the Euler-Lagrange equations? This leads us to the main results of this paper.

\begin{Thm}[Generic ill-posedness of the energy-momentum equations]
\label{thm:main}
Assume that $W\in C^\infty(GL(\R^n),\R)$. Furthermore, assume that $W(RX)=X$ for all $X\in \LL(\R^n)$ and $R\in O(\R^n)$. Then for every $\phi\in W^{1,\infty}(\Om, \R^n)$ such that $D\phi(x)\in \Int \Co O(\R^n)$ for almost every $x\in \Om$ the Dirichlet problem
\begin{equation}\label{eq:Thm1}
 \left\{
    \begin{array}{rl}
     \Div T(Du(x))=0 &\text{for a.e. }x\in \Om,\\
            u(x)=\phi(x) & x\in \dv \Om.
    \end{array} \right.
\end{equation}
has infinitely many solutions $u$ which can be taken to be nowhere $C^1$. 
\end{Thm}

This theorem shows that in order to hope for any type of partial regularity or uniqueness for solutions of the energy-momentum equations the assumption $\det(Du(x))>0$ a.e. or a smallness assumption on $\vert Du(x)\vert$ are essential. We can also ask if any solution \eqref{eq:Thm1} which is not $C^1$ can be a solution of the Euler-Lagrange equations? The next theorem shows that this cannot be the case under natural assumptions on the Lagrangian $W$. 

\begin{Thm}
\label{thm:main2}
Assume that $W\in C^1(\LL(\R^n))$ is frame indifferent (see Definition \ref{Def:FrameInd} below) and strictly rank-one convex. Then any weak solution $u\in W^{1,\infty}$ which is not $C^1$ to the differential inclusion $Du(x)\in \text{O}(\R^n)$ a.e. is not a weak solution of the Euler-Lagrange equations
\begin{align*}
\Div DW(Du(x))=0.
\end{align*}
\end{Thm}

\subsection*{ Previously known results}

The perhaps easiest non-smooth solution of the energy-momentum equation is the map $u:B_1(0)=\{x\in \R^3: \vert x\vert \leq 1\}\to S^2$ given by 
\begin{align*}
u(x)=\frac{x}{\vert x\vert}.
\end{align*}

This map belongs to $W^{1,2}(B_1(0),S^2)$ and is in fact (see \cite{L}) an absolute minimizer for the Dirichlet energy with its own boundary values in the space $W^{1,2}(B_1(0),S^2)$, a so called harmonic map to the unit sphere. As such it is also a weak solution of both the Euler-Lagrange equations and the energy-momentum equations. Due to the pointwise constraint $\vert u(x)\vert^2=1$, the Euler-Lagrange equations take the form 
\begin{align*}
\Delta u(x)=-\vert Du(x)\vert^2u(x).
\end{align*} 

However, the inner variations are fully compatible with pointwise constraints on the target, and are therefore the same as without the constraint. This is a general fact which can be used to construct various irregular solutions of the energy-momentum equations. Furthermore, since 

\begin{align*}
0=\nabla \vert u(x)\vert^2=2Du(x)u(x),
\end{align*} 
it follows that for all $x\neq 0$, $Du(x)$ is nowhere injective and $\det(Du(x))=0$ a.e.. Building on this example, and using a construction due to Ball and Murat in \cite{BM}, in \cite{SSp1} Sivaloganathan and Spector considered a class of frame indifferent (with respect to $SO(\R^n)$) and isotropic $W$ satisfying conditions that allow for cavitation solutions. They then construct a weak solution $u\in W^{1,p}(\R^n)$, $1\leq p<n$, of the energy-momentum equations such that $u$ has infinitely many discontinuities and yet is injective almost everywhere. They furthermore give an example of a $C^1$-solution of the energy-momentum equation which is not $C^2$. Finally, in \cite[Section 7]{Ti} Tione used convex integration theory to construct irregular solutions of the energy momentum equations of a specific functional. His method is however different from the present paper and give weaker results. 

In the positive direction and only restricted to dimension 2, in \cite{BOP} Bauman, Owen and Philips consider an energy density $W$ of the form 
\begin{align*}
W(X)=F(X)+H(\det(X)),
\end{align*} 
where $F$ is a \emph{quasi-convex} function and $H$ is a non-negative convex function. They show that any $C^{1,\alpha}$ solution of the energy-momentum equation in dimension 2 is in fact $C^{2,\alpha}$ and $\det(Du(x))>0$ for all $x$ in the domain.   

On the other hand also in dimension 2, Iwaniec, Kovalev and Oninnen prove in \cite{IKO} that for a large class of $W$, any solution $u\in W^{1,2}(\Om,\R^n)$ of the energy momentum equations which in addition satisfies $\det(Du(x))>0$ a.e. is in fact Lipschitz. See also \cite[Thm. 1.4 and Thm. 1.6]{MaYa} for cases when solutions of the energy momentum equations are in fact also homeomorphisms and unique. These regularity results for the energy-momentum equation are the strongest ones known to the author.

\subsection*{ Notation}

Let $\LL(\R^n,\R^m)$ denote the space of linear maps from $\R^n$ to $\R^m$. When $n=m$ we write $\LL(\R^n)$ instead of $\LL(\R^n,\R^n)$. For $X,Y\in \LL(\R^n,\R^m)$ we let $\vert X\vert=\sqrt{\text{tr}(X^\ast X)}$ denote the Hilbert-Schmidt norm, $\Vert X\Vert$ the operator norm and $\langle X,Y\rangle=\tr(X^\ast Y)$ the euclidean inner product. Furthermore, let
\begin{align*}
\text{GL}(\R^n)&=\{X\in \LL(\R^n): \det(X)\neq 0\},\\
\text{GL}_+(\R^n)&=\{X\in \LL(\R^n): \det(X)>0\},\\
\text{GL}_-(\R^n)&=\{X\in \LL(\R^n): \det(X)<0\},\\
O(\R^n)&=\{X\in \LL(\R^n): X^\ast X=I\},\\
SO_+(\R^n)&=\{X\in O(\R^n): \det(X)=1\},\\
SO_-(\R^n)&=\{X\in O(\R^n): \det(X)=-1\},\\
\text{Sym}_+(\R^n)&=\{X\in \LL(\R^n): X^\ast=X, X \geq 0\}. 
\end{align*}

If $X\in \LL(\R^n)$, we let $\text{cof}\,(X)$ denote the cofactor matrix of $X$ and $\text{adj}\,(X)=\text{cof}\,(X)^\ast$ denote the adjugate matrix of $X$. 

If $K\subset \R^n$ is a subset, then $\Co K$ denotes its convex hull, $\Rco K$ denotes its rank-convex hull and $\Int K$ denotes the interior of $K$.

\section*{ Energy-Momentum equations, symmetry and frame indifference}

\subsection*{ Energy-Momentum equations, symmetry and frame indifference}

\begin{Def}[Frame indifference]
\label{Def:FrameInd}
Let $W\in C^2(\LL(\R^n,\R^m))$ and consider the functional
\begin{align*}
\mathcal{W}[u]=\int_{\Om}W(Du(x))dx
\end{align*}
for an open set $\Om \subset \R^n$. We say that the functional is frame indifferent if $W(RX)=W(X)$ for all $R\in O(\R^m)$, where $O(\R^m)$ denotes the orthogonal group of $\R^m$.
\end{Def}

\begin{rem}
The reader should observe that the condition of frame indifference does not impose any material symmetry restrictions. Indeed frame indifference is just a manifestation of the fact that the energy should not change if an observer either rotates or reflects the coordinate system, i.e., if the mapping $u$ is changed to $Ru$ for some $R\in O(\R^n)$. Sometimes some authors only requires invariance under the special orthogonal group $SO(\R^n)$ rather than the full orthogonal group. However, it is natural to require invariance also under change of orientation of the coordinate system, at least for variational problems coming from classical physics. 
\end{rem}

\begin{Def}
The energy-momentum mapping $T: \LL(\R^n,\R^m)\to \LL(\R^n)$ associated to a Lagrangian $W\in C^\infty(\LL(\R^n,\R^m),\R)$ is defined according to
\begin{align}\label{eq:T2}
T(X):=X^\ast DW(X)-W(X)I_{\R^n},
\end{align}
where $I_{\R^n}$ is the identity mapping on $\R^n$. 
\end{Def}

When no confusion can arise we write $I$ instead of $I_{\R^n}$. 

\begin{Prop}[Energy-momentum mappings for frame indifferent Lagrangians]
\label{prop:Ind}
Let $W\in C^\infty(\LL(\R^n,\R^n),\R)$ be a frame indifferent Lagrangian, i.e., $W(RX)=W(X)$ for every $R\in O(\R^n)$. Then the energy-momentum mapping 
\begin{align*}
T(X)=X^\ast DW(X)-W(X)I
\end{align*}
is $O(\R^n)$-invariant, i.e. $T(RX)=T(X)$ for every $R\in O(\R^n)$. 
\end{Prop}

\begin{proof}
Let $X=RS$ be the left polar factorization of $X$, where $S=\sqrt{X^\ast X}$ and $R\in O(\R^n)$. Thus $W(X)=W(RS)=W(S)$. Thus there exists a function $\widetilde{W}(X^\ast X)=W(X)$ ($\widetilde{W}(X^\ast X):=W(\sqrt{X^\ast X})$). Then
\begin{align*}
W(X+\eps H)-W(X)&=\widetilde{W}((X+\eps H)^\ast (X+\eps H))-\widetilde{W}(X^\ast X)\\
&=\widetilde{W}(X^\ast X+\eps X^\ast H+\eps H^\ast X+\eps^2H^\ast H)-\widetilde{W}(X^\ast X)\\
&=\widetilde{W}(X^\ast X)+\eps\langle D\widetilde{W}(X^\ast X),X^\ast H+H^\ast X\rangle+o(\eps)-\widetilde{W}(X^\ast X)\\
&=\eps\langle D\widetilde{W}(X^\ast X),X^\ast H+H^\ast X\rangle+o(\eps). 
\end{align*} 

Furthermore using the cyclic invariance of the trace we find 
\begin{align*}
&\langle D\widetilde{W}(X^\ast X),X^\ast H+H^\ast X\rangle=\tr(D\widetilde{W}(X^\ast X)^\ast X^\ast H)+\tr(D\widetilde{W}(X^\ast X)^\ast H^\ast X)\\
&=\tr((XD\widetilde{W}(X^\ast X))^\ast H)+\tr(H^\ast XD_p\widetilde{W}(X^\ast X)^\ast)\\
&=\langle X(D\widetilde{W}(X^\ast X)+D\widetilde{W}(X^\ast X)^\ast),H\rangle.
\end{align*} 

Thus 
\begin{align*}
T(X)=X^\ast DW(X)-W(X)I=X^\ast X(D\widetilde{W}(X^\ast X)+D\widetilde{W}(X^\ast X)^\ast)-\widetilde{W}(X^\ast X)I.
\end{align*} 
This shows that $T(RX)=T(X)$ for all $R\in O(\R^n)$. 
\end{proof}

\begin{rem}\label{rem:DiffProp}
By \cite[Lemma 6.3, p.723]{Ball5} $\widetilde{W}$ is $C^\infty$ if $W$ is $C^\infty$ on the set $\{X\in \LL(\R^n): \det(X)\neq 0\}$. Therefore, in the proof of Proposition \label{prop:Ind} we have implicitly assumed that $\det(X)\neq 0$. The relation $T(RX)=T(X)$ for all $X$ is then deduced by continuity of $T$ and the fact that $\{X\in \LL(\R^n): \det(X)\neq 0\}$ is open and dense in $\LL(\R^n)$. 
\end{rem}

\begin{rem}\label{rem:Equi}
Note that $T(X)$ need not be a symmetric tensor. Furthermore, note that the structure field $ DW(X)=X(D\widetilde{W}(X^\ast X)+D\widetilde{W}(X^\ast X)^\ast)$ is {\bf not} $O(\R^n)-$invariant but $O(\R^n)$-equivariant, i.e. $DW(RX)=RDW(X)$ for all $R\in O(\R^n)$. Also note that in the case when $W$ is strictly convex, the map $DW:\LL(\R^n)\to \LL(\R^n)$ is strictly monotone and hence invertible. This is not case for $T$.
\end{rem}

\begin{Def}[Reduced energy-momentum tensor]
Let $W$ be a smooth frame indifferent Lagrangian and let $\widetilde{W}(X^\ast X)=W(X)$ for all $X\in \LL(\R^n)$. The \emph{reduced energy-momentum tensor} $\mathcal{T}$ is defined according to 
\begin{align*}
\mathcal{T}(X)=X^\ast X(D\widetilde{W}(X^\ast X)+D\widetilde{W}(X^\ast X)^\ast)-\widetilde{W}(X^\ast X)I.
\end{align*} 
\end{Def}
Note that since $X^\ast X$ is a symmetric non-negative linear map we think of the reduced energy-momentum tensor as defined only on the cone $\text{Sym}_+(n)\subset \LL(\R^n)$ of positive semidefinite linear maps.

\begin{Prop}[Symmetric energy-momentum tensors]
Let $W\in C^2(\LL(\R^n),\R)$ satisfy $W(XR)=W(X)$ for all $R\in SO(\R^n)$. Then $T(X)$ is symmetric, i.e, $T(X)^\ast=T(X)$.  
\end{Prop}

\begin{proof}
By the right polar factorization $X=SR$ where $S=\sqrt{XX^\ast}$ there exists  $\widetilde{W}$ such that $\widetilde{W}(XX^\ast)=W(X)$ for all $X$. Computing the gradient of $W$ using a similar computation as in the proof of Proposition \eqref{prop:Ind} shows that 
\begin{align*}
DW(X)=(D\widetilde{W}(X^\ast X)+D\widetilde{W}(X^\ast X)^\ast)X.
\end{align*} 
Thus, in this case the energy momentum tensor becomes
\begin{align*}
T(X)=X^\ast (D\widetilde{W}(X^\ast X)+D\widetilde{W}(X^\ast X)^\ast)X-\widetilde{W}(X)I
\end{align*}
which is symmetric. 
\end{proof}

\begin{rem}
Note that the case of symmetric energy-momentum tensor has been studied in the frame work of compensated integrability due to D. Serre in \cite{S} in the case hyperbolic evolution equations. 
\end{rem}

\section*{ Differential inclusions and generic ill-posedness}

\subsection*{ Differential inclusions in the orthogonal group}

In this section we will consider some background material concerning differential inclusions into the orthgonal group. Consider the relaxation of the orthogonal group $O(\R^n)$, i.e., its convex hull given by 
\begin{align*}
\text{co}\,O(\R^n)=\{X\in \LL(\R^n): X^\ast X\leq I\}. 
\end{align*}

\begin{Thm}\label{thm:OrtDiff}
Let $\Om\subset \R^n$ be a Lipschitz domain and let $\phi\in W^{1,\infty}(\Om,\R^n)$ be a Lipschitz map such that $D\phi(x)\in \Int \Co O(\R^n)$ for a.e. $x\in \Om$. Then the differential inclusion 
\begin{align}\label{eq:DiffO}
 \left\{
    \begin{array}{rl}
      Du(x)\in O(\R^n) & \text{for a.e. $x\in \Om$ },\\
      u(x)=\phi(x)& \text{for a.e. $x\in \dv\Om$ }
     \end{array} \right.
\end{align}
for $u\in W^{1,\infty}(\Om,\R^n)$ possesses infinitely many solutions that are nowhere $C^1$. 
\end{Thm}

The key to this result is that $O(\R^n)$ posses many rank one-connections, i.e. there exits matrices $X,Y\in O(\R^n)$ such that 
\begin{align*}
X-Y=u\otimes v
\end{align*}
for some $u,v\in \R^n$. 

A proof of Theorem \ref{thm:OrtDiff} can be found in \cite[Sec. 5.1]{Sz2}.

Solutions of differential inclusions need not necessarily be very irregular. 
\begin{Def}\label{Def:Part}
Let $u\in C^{0,1}(\Om,\R^m)$ be a Lipschitz map. Let 
\begin{align}
\Sigma(u)=\{x\in \Om: \text{$u$ is not $C^1$ in a neighbourhood of $x$}\}
\end{align}
denote the \emph{singular set} of $u$. 
\end{Def}

In \cite{DaMP}, the authors consider solutions of \eqref{eq:DiffO} with affine boundary values generated by orgami maps. These maps are piecewise $C^1$ and the Hausdorff measure $\mathscr{H}^{n-1}(\Sigma(u))$ is locally finite in the interior of the domain. However to satisfy the boundary conditions $\Sigma(u)$ will become fractal like as we approach the boundary and $\mathscr{H}^{n-1}(\Sigma(u))=+\infty$ in the whole domain. Similar types of solutions are also considered in the paper \cite{IVV} and also for the energy-momentum equations for the Dirichlet energy in \cite[Sec. 3.6]{IO}. In both cases $Du(x)\in K$, where $K$ is a \emph{finite subset} of $O(\R^n)$.

\subsection*{ Generic ill-posedness}

As we have seen in Proposition \ref{prop:Ind} the mapping $T: \LL(\R^n)\to \LL(\R^n)$ is not invertible for a frame indifferent Lagrangian $W$. In particular the level sets of
\begin{align*}
T(X)=Y
\end{align*}
for some fixed $Y\in T(\LL(\R^n))$ are $O(\R^n)$-invariant. 

\begin{Lem}\label{lem:key}
Any solution $Du(x)\in O(\R^n)A$ for some $A\in \LL(\R^n)$ is a solution of the energy-momentum equations.
\end{Lem}

\begin{proof}

If $Du(x)\in O(\R^n)A$ a.e. then 
\begin{align*}
T(Du(x))=T(A)=Y
\end{align*}
is constant a.e. and hence a weak solution of $\text{div}\,T(Du(x))=0$. 
\end{proof}

This leads to differential inclusions of the form 
\begin{align*}
Du(x)\in T^{-1}(Y)
\end{align*}
for some fixed $Y$ such that $T^{-1}(Y)\neq \varnothing$.

\begin{proof}[Proof of Theorem \ref{thm:main}]
The proof follows by combining Lemma \ref{lem:key} with Theorem \ref{thm:OrtDiff}.
\end{proof}

There are similar results for the lack partial regularity for the Euler-Lagrange equations for elliptic systems. In \cite[Theorem 4.1]{MS2} the authors show that there exists a smooth strongly quasiconvex function $W:\LL(\R^2)\to \R$ such that there exists a Lipschitz continuous solution of $\text{div}DW(Du(x))=0$ that is nowhere $C^1$. This is done by rewriting the Euler-Lagrange equation as a differential inclusion and using methods from convex integration theory. The result was extended in \cite{Sz} to apply also to smooth strongly polyconvex functions. Moreover, the solutions in \cite{Sz} are also weak local minimisers. This also applies to the example in \cite{MS2} by the work of \cite{KT}. It is however important to note that weak local minimisers $u_0\in W^{1,\infty}$ need not be weak solutions of the energy-momentum equations. Indeed, let $\Om\subset \R^n$ be a domain and let $\lambda\in C^\infty_0(\Om,\R^n)$. Let $u_0\in W^{1,\infty}(\Om,\R^n)$ and let $u_\eps(x)=u_0(x+\eps \lambda(x))$ be an inner variation. Then 
\begin{align*}
\vert Du_0(x)-Du_\eps(x)\vert&=\vert Du_0(x)-Du_0(x+\eps \lambda(x))(I+\eps D\lambda (x))\vert\\&\geq \big\vert \vert Du_0(x)-Du_0(x+\eps \lambda(x))\vert -\eps \vert Du(x+\eps\lambda(x))D\lambda (x)\vert\big\vert\\&\geq \big\vert \vert Du_0(x)-Du_0(x+\eps \lambda(x))\vert -\eps \Vert \vert Du\vert \Vert_{L^\infty(\Om)}\Vert  \vert D\lambda\vert \Vert_{L^\infty(\Om)}\big\vert
\end{align*} 
Since $Du_0$ is not continuous it may happen that \newline$\vert Du_0(x)-Du_0(x+\eps \lambda(x))\vert\geq 1+\eps \Vert \vert Du\vert\Vert_{L^\infty(\Om)}\Vert \vert D\lambda\vert \Vert_{L^\infty(\Om)}$ for every $\eps>0$ and thus that $\Vert u_0-u_\eps\Vert_{W^{1,\infty}(\Om)}>1$ for all $\eps>0$. Thus inner variations need not be close in $W^{1,\infty}(\Om)$-norm. Therefore $u_0$ being a weak local minimiser need not imply that $u_0$ solves the energy-momentum equations.

It is therefore a natural question if the differential inclusions giving solutions to the energy-momentum equations are also weak solutions of the Euler-Lagrange equations? For this purpose we first consider a special class of solutions given by laminations, the reason being that any solution in for example \cite{DaMP,IVV} is locally a lamination outside a small closed set. Given a first order partial differential operator $\mathscr{A}$ with constant coefficients and its associated symbol $\mathbb{A}(\xi)$ with $\xi \in \mathbb{S}^{n-1}$ consider functions of the form 
\begin{align}\label{eq:lam}
u(x)=\lambda h(\langle x,\xi\rangle)+\mu(1-h(\langle x,\xi\rangle))
\end{align}
where $h:\R\to \{0,1\}$ is measurable and $\lambda-\mu\in \text{ker}\,\mathbb{A}(\xi)$. These are solutions to two state rigidity problem

\begin{equation*}
 \left\{
    \begin{array}{ll}
     \mathscr{A}u(x)=0 &\text{in the sense of distributions},\\
            u(x)\in \{\mu,\lambda\}.&
    \end{array} \right.
\end{equation*}

In the case when $\mu,\lambda\in O(\R^n)$ are rank-1 connected then the laminate solution \eqref{eq:lam} gives a solution to the differential inclusion $Du(x)\in O(\R^n)$ a.e. with $\mathcal{A}=\text{curl}$, where  $\text{curl}$ denotes the matrix curl operator. In this case the boundary values for non-trivial measurable functions $h$ are however, not smooth. Furthermore, we want the rank-1 connected laminate to be such that $v(x)=DW(D(x))$ is two state laminate solution for $\mathscr{A}=\text{div}$, where $\text{div}$ is the matrix divergence. By \cite[p. 8]{KR} $\mathbb{A}(\xi)X=X\xi$ for an $n\times n$-matrix $X$ and hence we must have that $(DW(\mu)-DW(\lambda))\xi=0$ as well. This leads us to the following proposition.

\begin{Prop}\label{prop:Lam}
Let $W$ satisfy the assumptions in Theorem \ref{thm:main}. Let $A,B\in O(\R^n)$ be rank one connected and consider the laminate 
\begin{align}\label{eq:Lamsol}
Du(x)=A h(\langle x,\xi\rangle)+B(1-h(\langle x,\xi\rangle))
\end{align}
for some measurable $h:\R\to \{0,1\}$ and such that $B-A=a\otimes \xi$ for some $a\in \R^n$ and $\xi \in \mathbb{S}^{n-1}$.
Then the laminate is a distributional solution to the Euler-Lagrange equation
\begin{align*}
\Div DW(Du(x))=0
\end{align*}
if and only if $\langle \xi, T(I)\xi\rangle=-W(I)$ or equivalently if and only if $\langle \xi, DW(I)\xi\rangle=0$. 
\end{Prop}

\begin{proof}
Let $\lambda=DW(A)$ and $\mu=DW(B)$. Then for the matrix field $M(x)=DW(Du(x))$ to be divergence free in the sense of distributions we must have $(DW(B)-DW(A))n=0$. Since 
$O(\R^n)$ acts transitively on itself there exists a $U\in O(\R^n)$ such that 
\begin{align*}
B=UA.
\end{align*}
On the other hand we have 
\begin{align*}
B-A=a\otimes \xi
\end{align*}
which implies 
\begin{align*}
(U-I)A=a\otimes \xi \quad \Longrightarrow \quad U-I=  a\otimes \xi\circ A^\ast
\end{align*}
Using that $DW$ is $O(\R^n)$-equivariant by Remark \ref{prop:Ind}, we find that 
\begin{align*}
(DW(B)-DW(A))\xi&=(UDW(A)-DW(A))\xi=(U-I)DW(A)\xi\\
&= a\otimes \xi\circ A^\ast\circ  DW(A)\xi =a\langle \xi, A^\ast DW(A)\xi\rangle\\
&=a\langle \xi, (T(A)+W(A)I)\xi\rangle=a(\langle \xi, T(A)\xi\rangle+W(A)\vert \xi\vert^2)\\&=a(\langle \xi, T(A)\xi\rangle+W(A))
\end{align*}
Thus $\langle \xi, T(A)\xi\rangle+W(A)=0$. Since $A\in O(\R^n)$ and $T$ and $W$ are $O(\R^n)$-invariant $\langle \xi, T(A)\xi\rangle+W(A)=\langle \xi, T(I)\xi\rangle+W(I)=0$. 
\end{proof}

\begin{rem}\label{rem:Ela}
If $W$ is the energy density of a hyperelastic material it is physically reasonable that scalings $x\mapsto tx$ for $t>0$ costs energy. Consequently the function 
\begin{align*}
j(t)=W(Du(x))=W(tI)
\end{align*}
should have a minimum at $t=1$. Since
\begin{align*}
j'(t)=\langle DW(tI),tI\rangle
\end{align*}
we find that the condition $j'(1)=0$ implies that $0=\tr(DW(I))=\tr(T(I)+W(I)I)=\tr(T(I))+nW(I)$. Since we have for an ON-basis $\{e_j\}_{j=1}^n$
\begin{align*}
\tr(T(I))+nW(I)=\sum_{j=1}^n(\langle e_j,T(I)e_j\rangle +W(I)),
\end{align*}
we see that the condition $j'(1)=0$ can be seen as an averaged condition of the previous Proposition \ref{prop:Lam}.
\end{rem}

It is easy to produce frame indifferent smooth functions $W$ which satisfy $DW(I)=0$, $W(X)=(\vert X\vert^2-4)^2$ for example will do. However $W$ is not rank-one convex and therefore not polyconvex either. In Example \ref{ex:Qmean} we give a less trivial example coming from geometric function theory and in particular the study of mappings of finite distortion. This functional is not globally polyconvex but polyconvex when restricted to $GL_+(\R^n)$.

\begin{Prop}
Assume that $W\in C^1(\LL(\R^n))$ is frame indifferent and strictly rank-one convex. Let $u$ be a laminate solution to the differential inclusion $Du(x)\in \text{O}(\R^n)$ as in Proposition \ref{prop:Lam}. Then $u$
is not a weak solution to the Euler-Lagrange equations.
\end{Prop}

\begin{proof}
By a previous remark it is sufficient to consider solutions with $Du(x)\in \{I,R\}$ and $R\in O(\R^n)$ rank-one connected to $I$. Let $R-I=a\otimes \xi$ for some $\xi \in S^{n-1}$ and some $a\in \R^n$. Furthermore since $R^\ast R=I$ we get the equation 
\begin{align*}
\xi \otimes a+a\otimes \xi=-\vert a\vert^2\xi\otimes \xi
\end{align*}
which implies that $a=-2\xi$. 

For $u$ to solve the Euler-Lagrange equations we must have using the $\text{O}(\R^n)$-equivariance of $DW$ 
\begin{align*}
0=(DW(R)-DW(I))\xi=(R-I)DW(I)\xi=2(\xi\otimes \xi)\circ (DW(I)\xi)=2\xi\langle \xi,DW(I)\xi\rangle
\end{align*}
and so $\langle \xi,DW(I)\xi\rangle=0$. 
Now consider the function 
\begin{align*}
\phi(t)=W(I-2t\xi \otimes \xi)
\end{align*}
By frame indifference $\phi(0)=\phi(1)$ and
\begin{align*}
\phi'(t)=\langle DW(I-2t\xi \otimes \xi),R-I\rangle=-2\langle DW(I-2t\xi \otimes \xi),\xi \otimes \xi \rangle.
\end{align*}
By assumption $\phi(t)$ is strictly convex. Hence $\phi'(0)=\phi'(1)\neq 0$. Thus 
\begin{align*}
0\neq \langle DW(I),\xi \otimes \xi \rangle=\langle \xi,DW(I)\xi\rangle
\end{align*}
a contradiction. 
\end{proof}

\begin{proof}[Proof of Theorem \ref{thm:main2}]
We note that by the frame indifference of $W$ (Remark \ref{rem:Equi}) $DW(I)^\ast=DW(I)$. Furthermore since any rank-one connected matrix $R\in \text{O}(\R^n)$ to $I$ is of the form $I-2\xi\otimes \xi$ for some $\xi\in S^{n-1}$ we have for $\phi_\xi(t)=W(I-2t\xi \otimes \xi)$
\begin{align*}
\phi_\xi'(t)=-2\langle DW(I-2t\xi \otimes \xi), \xi\otimes \xi\rangle
\end{align*}
and the strict rank-one convexity assumption implies that $\phi_\xi'(0)=\phi'_\xi(1)\neq 0$ for all $\xi \in S^{n-1}$ we find that the quadratic form $Q(\xi)=\langle DW(I), \xi\otimes \xi\rangle=\langle \xi,DW(I)\xi \rangle\neq 0$
for all $\xi\in S^{n-1}$. Thus either $Q(\xi)>0$ or $Q(\xi)<0$. We may assume the former case. Hence $DW(I)$ is positive definite. Set $A=DW(I)$ and define the linear map $L\in \LL(\LL(\R^n))$ by 
\begin{align*}
L(X)=XA.
\end{align*}
Then $L$ is symmetric and positive definite. Indeed,
\begin{align*}
\langle Y,L(X)\rangle&=\tr(Y^\ast XA)=\tr(AY^\ast X)=\tr((YA^\ast)^\ast X)=\tr((YA)^\ast X)\\&=\langle L(Y),X\rangle.
\end{align*}
Furthermore, using that $A$ is diagonalizable with diagonal matrix $D$ and eigenvalues $0<\lambda_1\leq \lambda_2\leq ...\leq \lambda_n$ such that $A=RDR^\ast$ and setting $Y=XR$ we find 
\begin{align*}
\langle X,L(X)\rangle&=\tr(X^\ast XRDR^\ast)=\tr(R^\ast X^\ast XRD)=\tr((X^\ast R)^\ast XRD)\\
&=\langle Y,YD\rangle=\sum_{j=1}^n\lambda_j \langle Y,Y(e_j\otimes e_j)\rangle=\sum_{j=1}^n\lambda_j \langle Y(e_j),Y(e_j)\rangle\\
&\geq \lambda_1 \sum_{j=1}^n\langle Y(e_j),Y(e_j)\rangle=\lambda_1 \vert Y\vert^2=\lambda_1\vert X\vert^2. 
\end{align*}
Now assume that $u$ is a solution of the differential inclusion $Du(x)\in O(\R^n)$ that is not $C^1$ and that in addition $u$ is a weak solution of the Euler-Lagrange equations on some domain $\Om \subset \R^n$. Then for any $\varphi\in C^\infty_0(\Om, \R^n)$ and using the $\text{O}(\R^n)$-equivariance of $W$ we find 
\begin{align*}
0=\int_{\Om}\langle DW(Du(x)),D\varphi(x)\rangle dx=\int_{\Om}\langle Du(x)DW(I),D\varphi(x)\rangle dx=\int_{\Om}\langle L(Du(x)),D\varphi(x)\rangle dx.
\end{align*}
Thus $u$ is a weak solution of the very strongly elliptic constant coefficient equation (in the sense of \cite[Definition 3.36 (3.16), p.53]{GM})
\begin{align*}
\text{div}\,L(Du(x))=0.
\end{align*}
However, by elliptic regularity theory (\cite[Thm. 4.11]{GM}) $u\in C^\infty$, a contradiction. 
\end{proof}

\begin{rem}
Note that the assumption that $W$ is $C^1$ and strictly rank-one convex implies that $DW(I)$ is positive definite is incompatible with $\tr(DW(I))=0$, which is a natural condition in nonlinear elasticity. See the discussion in Remark \ref{rem:Ela}. 
\end{rem}

\begin{rem}
There has been other attempts in \cite{DeDeKT} to construct stationary points of strictly polyconvex functionals by extending the methods in \cite{MS2}. The main result of \cite{DeDeKT} is that the methods do not extend to this case, giving further indication that stationary points may in fact possess some form of partial regularity. 
\end{rem}

\section*{ Invertibility of the reduced energy-momentum mapping and double well inclusion}

So far we have considered solutions $u\in W^{1,\infty}$ to the differential inclusion $Du(x)\in O(\R^n)$ a.e. All these solutions have in common that the essential range of $\det(Du(\Om))$ lies in $\{-1,1\}$. One can ask, in particular with respect to the results in \cite{MS1}, whether it is possible to find other types of differential inclusions which are also solutions of the energy-momentum equations and such that $\det(Du(x))>0$. Moreover one can ask if the Lipschitz regularity result in \cite{IKO} can be improved to a partial regularity result, i.e, if one can show that the singular set $\Sigma_u$ of Definition \ref{Def:Part} is a closed set with $\mathscr{H}^{n-1}(\Sigma_u)=0$. We now formulate an obstruction to such a result.

Indeed, if we assume that the reduced energy momentum tensor $\mathcal{T}:\text{Sym}_+(n)\to \LL(\R^n)$ is not injective then we could find two solutions $A,B\in \text{Sym}_+(n)$ such that 
\begin{align*}
\mathcal{T}(A)=\mathcal{T}(B)=S
\end{align*}
with $\det(A)>0$ and $\det(B)>0$. In particular we would have 
\begin{align*}
T(SO(n)A)=T(SO(n)B)=S. 
\end{align*}
and solutions of the differential inclusion $Du(x))=SO(n)A\cup SO(n)B$ for a.e. $x\in \Om$ would also be solutions of the energy-momentum equations. This differential inclusion is studied in \cite{MS1,Da} in the case when $n=2$.  The following theorem holds.
\begin{Thm}\cite[Cor. 1.4 ]{MS1}, \cite[Thm. 10.28]{Da}
\label{thm:DW}
Let $A,B\in \LL(\R^2)$ be diagonal matrices whose diagonal entries are $a_1,a_2$ and $b_1,b_2$ respectively. Assume that $0<b_1< a_1\leq a_2<b_2$ and $\det(A)\leq \det(B)$. Let $T\in \Int \Rco (SO_+(\R^2)A\cup SO_+(\R^2)B)$. Let $\phi(x)=Tx+c$, with $c\in \R^n$. Then the differential inclusion
\begin{equation}\label{eq:DiffDet}
 \left\{
    \begin{array}{lr}
     Du(x))\in SO_+(\R^2)A\cup SO_+(\R^2)B& \text{for a.e. }x\in \Om,\\
            u(x)=\phi(x) & x\in \dv \Om.
    \end{array} \right.
\end{equation}
has infinitely many solutions. 
\end{Thm}

\begin{rem}
In Theorem \ref{thm:DW} the assumption that $A$ and $B$ are diagonal matrices is not essential, as one can always reduce to this case. Also note the importance that $A\neq \alpha I$ and $B\neq \beta I$ for some constants $\alpha>0,\beta>0$, since \eqref{eq:DiffDet} in this case implies that $u$ is a conformal map and by Liouville's theorem \cite[Thm. 5.1.1]{IM} any conformal map $u\in W^{1,n}$ is a Möbius transformation of $\R^n\cup\{\infty\}$. 
\end{rem}

\begin{rem}
Note that solutions of \eqref{eq:DiffDet} need in no way be locally injective. In particular since both $I$ and $-I$ belongs to $\text{SO}_+(\R^2)$ and $0=\frac{1}{2}I-\frac{1}{2}I\in \text{co}(\text{SO}_+(\R^2))$ the inverse function theorem for Lipschitz mappings due to Clarke (\cite{Clarke}) does not apply. 
\end{rem}

If a convex, polyconvex or quasiconvex $W$ whose reduced energy-momentum tensor is not an injective map, and such that there exists diagonal matrices $A,B\in \text{Sym}_+(\R^n)$ that satisfies the assumption of Theorem \ref{thm:DW}, then Theorem \ref{thm:DW} would show (at least in dimension two) that there are energy momentum tensors for which well-posedness of the Dirichlet problem for the energy-momentum equations fails, even with the additional constraint $\det(Du(x))>0$ for a.e. $x$ and furthermore no partial $C^1$-regularity holds. In all the examples we study there are cases when the reduced energy momentum tensor fails to be injective. The solutions to $\mathcal{T}^{-1}(Z)$ however fails to satisfy the conditions of Theorem \ref{thm:DW}. 
\begin{OP}
Does there exists a smooth strictly convex, polyconvex or quasiconvex $W\in C^\infty(\LL(\R^2),\R)$ for which its reduced energy-momentum tensor fails to be an injective map and for which there exists matrices $A$ and $B$ that satisfy the assumptions of Theorem \ref{thm:DW} and such that $\mathcal{T}(A)=\mathcal{T}(B)$?
\end{OP}

\subsection*{ Invertibility of the reduced energy-momentum tensor in a number of interesting cases}

In this section we will consider a number of important functionals that occur in nonlinear elasticity and geometric functions theory. We will show that in all these cases the reduced energy-momentum mapping $\mathcal{T}$ is typically not injective, yet the conditions of Theorem \ref{thm:DW} are not satisfied. In addition, they all have the feature that their Lagrangian $W$ in addition to being frame indifferent is also isotropic, i.e., 
\begin{align*}
W(R^\ast XR)=W(X)
\end{align*}
 for all $R\in O(\R^n)$ and all $X\in \LL(\R^n)$. We begin with the Dirichlet $p$-energy.

\begin{ex}
\label{ex:Convex}
Let $1<p<+\infty$ and let for $u\in W^{1,p}(\Om,\R^n)$
\begin{align*}
\mathscr{D}_p[u]=\int_{\Om}\vert Du(x)\vert^p dx.
\end{align*}

Since $W(X)=\vert X\vert^p$ we find that $DW(X)=p\vert X\vert^{p-2}X$ and the energy-momentum mapping becomes 
\begin{align*}
T(X)=p\vert X\vert^{p-2}X^\ast X-\vert X\vert^pI=p\text{tr}(X^\ast X)^{(p-2)/2}X^\ast X-\text{tr}(X^\ast X)^{p/2}I. 
\end{align*}
The reduced energy-momentum mapping becomes with $Y=X^\ast X$
\begin{align*}
\mathcal{T}(Y)=p\text{tr}(Y)^{(p-2)/2}Y-\text{tr}(Y)^{p/2}I. 
\end{align*}

Note that 
\begin{align*}
\text{tr}(\mathcal{T}(Y))=p\text{tr}(Y)^{(p-2)/2}\text{tr}(Y)-\text{tr}(Y)^{p/2}\text{tr}(I)=(p-n)\text{tr}(Y)^{p/2},
\end{align*}
which is different from $0$ if and only if $p\neq n$. Let $Z\in \mathcal{T}(\text{Sym}_+(n))$. Consider the equation for $p\neq n$
\begin{align*}
p\text{tr}(Y)^{(p-2)/2}Y-\text{tr}(Y)^{p/2}I=Z. 
\end{align*}
Taking traces of both sides gives us 
\begin{align*}
(p-n)\text{tr}(Y)^{p/2}=\text{tr}(Z)\quad \Longrightarrow \quad \text{tr}(Y)=((p-n)^{-1}\text{tr}(Z))^{2/p}.
\end{align*}
Thus
\begin{align}\label{eq:Trace}
Y=\frac{Z+\text{tr}(Y)^{p/2}I}{p\text{tr}(Y)^{p/2-1}}=\frac{Z+\frac{1}{p-n}\text{tr}(Z)I}{p(\frac{1}{p-n}\text{tr}(Z))^{1-2/p}}. 
\end{align}
Since $Y$ is positive semi-definite it has a unique positive semi-definite square root $\sqrt{Y}$. In particular, all solutions of $T(X)=Z$ are given by 
\begin{align*}
T^{-1}(Z)=O(\R^n)\sqrt{\frac{Z+\frac{1}{p-n}\text{tr}(Z)I}{p(\frac{1}{p-n}\text{tr}(Z))^{1-2/p}}}
\end{align*}
and the situation \eqref{eq:DiffDet} cannot occur. Furthermore, by Theorem \ref{thm:main2} and in view of Uhlenbeck's regularity result \cite{U}, weak solutions $u\in W^{1,p}(\Om,\R^n)$ of the Euler-Lagrange equations of the Dirichlet $p$-energy are always $C^{1,\alpha}$ for some $0<\alpha<1$. Hence the weak solutions of the energy-momentum equations which are not $C^{1,\alpha}$ are not weak solutions of the Euler-Lagrange equations. In the conformally invariant case $p=n$ the formula \eqref{eq:Trace} does not hold and in fact we now show that $\mathcal{T}(Y)$ is not injective. If follows from the equation $\mathcal{T}(Y)=Z$ that $[Y,Z]=0$ so if $Z$ is diagonalizable so is $Y$. Thus we restrict to considering only diagonal matrices $Y$ and only consider the case $n=2$. We find that if 
\begin{align*}
Y=\begin{bmatrix}
\alpha & 0\\
0 & \beta
\end{bmatrix},\quad
Z=\begin{bmatrix}
c & 0\\
0 & -c
\end{bmatrix}
\end{align*}
where $c\geq 0$ we find the system of equations 
\begin{equation*}
 \left\{
    \begin{array}{l}
\alpha-\beta=c \\
\beta-\alpha=-c
    \end{array} \right.
\end{equation*}
Hence if $\alpha=t$, $\beta=t-c$ and $t\geq c$ parametrises the solutions. Thus $\mathcal{T}$ is not an injective map. On the other hand we can find no two $t_1,t_2$ such that the condition of Theorem \ref{thm:DW} is satisfied. 
\end{ex}

\begin{ex}[ $q$-mean distortion]
\label{ex:Qmean}

Let $u:\Om \subset \R^n \to \R^n$ be a map in $W^{1,n}(\Om,\R^n)$ and consider the $q$-mean distortion functional

\begin{align*}
\mathscr{K}_q[u]=\int_{\Om}\mathbb{K}(x,u)^qdx=\int_{\Om}\bigg(\frac{\vert Du(x)\vert^n}{J(x,u)}\bigg)^qdx
\end{align*}
where $q\geq 1$. $W$ is a priori only well-defined when $\det(X)>0$, however we can extend $W$ to $\widetilde{W}\in C^\infty(\text{GL}_+(\R^n)\cup \text{GL}_-(\R^n))$ as a frame indifferent function by defining 
\begin{align*}
\widetilde{W}(X)=\bigg(\frac{\vert X\vert^n}{\vert \det(X)\vert}\bigg)^q. 
\end{align*}
This extension is however not polyconvex due to the blow up when $\det(X)=0$ (except in the case $X=tI$ and $t\to 0$). The only polyconvex extension is to define $\mathbb{K}(x,u)^q=+\infty$ whenever $\det(X)\leq 0$. For more on this functional we refer the reader to \cite{IMO} and references there in. We will however use $\widetilde{W}$ and by abuse of notation also write $W$ for its frame indifferent extension. Set $W(X)=(f(X)g(X))^q$ where 
$f(X)=\vert X\vert^n$ and $g(X)=\vert \det(X)\vert^{-1}$. Then 
\begin{align*}
DW(X)=q\bigg(\frac{\vert X\vert^n}{\vert \det(X)\vert}\bigg)^{q-1}[g(X)Df(X)+f(X)Dg(X)]. 
\end{align*}

By \eqref{eq:Det} 
\begin{align*}
Df(X)&=n\vert X\vert^{n-2}X,\\
Dg(X)&=-\frac{\text{sgn}(\det(X))}{\vert \det(X)\vert^2}\text{adj}(X)^\ast
\end{align*}

Thus,
\begin{align*}
DW(X)=q\bigg(\frac{\vert X\vert^n}{\vert \det(X)\vert}\bigg)^{q-1}\bigg[\frac{n\vert X\vert^{n-2}X}{\vert \det(X)\vert}-\frac{\text{sgn}(\det(X))\vert X\vert^n}{\vert \det(X)\vert^2}\text{adj}(X)^\ast\bigg]
\end{align*}
and
\begin{align*}
T(X)&=X^\ast DW(X)-W(X)I\\
&=q\bigg(\frac{\vert X\vert^n}{\vert \det(X)\vert}\bigg)^{q-1}\bigg[\frac{n\vert X\vert^{n-2}X^\ast X}{\vert \det(X)\vert}-\frac{\text{sgn}(\det(X))\vert X\vert^n}{\vert \det(X)\vert^2}X^\ast \text{adj}(X)^\ast\bigg]-\bigg(\frac{\vert X\vert^n}{\vert \det(X)\vert}\bigg)^qI\\
&=q\bigg(\frac{\vert X\vert^n}{\vert \det(X)\vert}\bigg)^{q-1}\bigg[\frac{n\vert X\vert^{n-2}X^\ast X}{\vert \det(X)\vert}-\frac{\vert X\vert^n}{\vert \det(X)\vert}I\bigg]-\bigg(\frac{\vert X\vert^n}{\vert \det(X)\vert}\bigg)^qI
\end{align*}

Thus,

\begin{align*}
T(I)&=q\bigg(\frac{\vert I\vert^n}{\vert \det(I)\vert}\bigg)^{q-1}\bigg[\frac{n\vert I\vert^{n-2}I}{\vert \det(I)\vert}-\frac{\vert I\vert^n}{\vert \det(I)\vert}I\bigg]-\bigg(\frac{\vert I\vert^n}{\vert \det(I)\vert}\bigg)^qI\\
&=-W(I)I
\end{align*}
and $\langle T(I)\xi,\xi\rangle=-W(I)$. 

Hence $T$ does satisfy the assumptions of Proposition \ref{prop:Lam} and there are stationary points of the functional which are nowhere $C^1$. Indeed as an explicit example take $n=2$ and $q=1$ and let 
\begin{align*}
R=\begin{bmatrix}
1 & 0\\
0 &-1
\end{bmatrix}.
\end{align*}
Then $I-R=2e_2\otimes e_2$ and so $R$ is rank-one connected to $I$. Furthermore, 
\begin{align*}
DW(I)=\frac{2I}{\vert \det(I)\vert}-\frac{\text{sgn}(\det(I))\vert I\vert^2}{\vert \det(I)\vert^2}\text{adj}(I)^\ast=0,
\end{align*}
and
\begin{align*}
\text{adj}(R)^\ast=\begin{bmatrix}
-1 & 0\\
0 &1
\end{bmatrix},
\end{align*}
and 
\begin{align*}
DW(R)=\frac{2R}{\vert \det(R)\vert}-\frac{\text{sgn}(\det(R))\vert R\vert^2}{\vert \det(R)\vert^2}\text{adj}(R)^\ast=2[R+\text{adj}(R)^\ast]=0.
\end{align*}
Thus, $u$ is also a solution to the Euler-Lagrange equations. 
The reduced energy-momentum tensor becomes
\begin{align*}
\mathcal{T}(Y)&=q\bigg(\frac{\text{tr}(Y)^{n/2}}{ \sqrt{\det(Y)}}\bigg)^{q-1}\bigg[\frac{n(\text{tr}(Y))^{(n-2)/2}Y}{\sqrt{\det(Y)}}-\frac{\text{tr}(Y)^{n/2}}{\sqrt{\det(Y)}}I\bigg]-\bigg(\frac{\text{tr}(Y)^{n/2}}{\sqrt{\det(Y)}}\bigg)^qI\\
&=q\bigg(\frac{\text{tr}(Y)^{n/2}}{ \sqrt{\det(Y)}}\bigg)^{q}\bigg[\frac{nY}{\tr(Y)}-I\bigg]-\bigg(\frac{\text{tr}(Y)^{n/2}}{\sqrt{\det(Y)}}\bigg)^qI\\
&=\bigg(\frac{\text{tr}(Y)^{n/2}}{ \sqrt{\det(Y)}}\bigg)^{q}\bigg[\frac{qnY-(q+1)\text{tr}(Y)I}{\tr(Y)}\bigg],
\end{align*}
provided $\det(Y)\neq 0$ and $\tr(Y)\neq 0$. We now specialise to the case when $n=2$ and $q=1$. Then
\begin{align*}
\text{tr}(\mathcal{T}(Y))&=-2\bigg(\frac{\text{tr}(Y)}{\sqrt{\det(Y)}}\bigg)\\
\det(\mathcal{T}(Y))&=\det\bigg[\frac{2}{\sqrt{\det(Y)}}(Y-\tr(Y)I)\bigg]=\frac{4}{\det(Y)}\det(Y-\tr(Y)I)\\
&=\frac{2}{\det(Y)}((\tr(Y-\tr(Y)I))^2-\tr((Y-\tr(Y)I)^2))\\
&=\frac{2}{\det(Y)}((\tr(Y))^2-\tr(Y^2-2\tr(Y)Y+(\tr(Y))^2I))\\
&=\frac{2}{\det(Y)}((\tr(Y))^2-\tr(Y^2))=\frac{4}{\det(Y)}\det(Y)=4
\end{align*}
Assume that $Z\in \text{im}(\mathcal{T})$ and that furthermore $Y$ is a diagonal matrix. Then we recall that $Y-\tr(Y)I=-\adj(Y)$. This gives us the equation
\begin{align*}
\frac{2}{\sqrt{\det(Y)}}[Y-\tr(Y)I]=Z,\,\,\, \Longleftrightarrow -2\frac{\adj(Y)}{\sqrt{\det(Y)}}=Z
\end{align*}
or upon multiplying by $Y$,
\begin{align}\label{eq:Ysolv}
-2\sqrt{\det(Y)}I=YZ. 
\end{align}
Set $t=\sqrt{\det(Y)}$. This gives the one-parameter family of solutions
\begin{align}\label{eq:Ysolv2}
Y(t)=-2tZ^{-1}
\end{align}
for $t>0$ given a solution $Y_0$ such that $\mathcal{T}(Y_0)=Z$. For example if $Y_0=I$, then $Z=-2I$, and $Y(t)=tI$.  
Again one may verify that no two solutions of \eqref{eq:Ysolv2} satisfy the assumptions of Theorem \ref{thm:DW}. 
\end{ex}

\begin{ex}
\label{ex:NonInjec}

Let $W:\LL(\R^n)\to \R$ be given by 
\begin{align*}
W(X)&=\vert X\vert^p+\vert X^{-1}\vert^p\vert \det(X)\vert.
\end{align*}
This is an example from \cite{IMO} and comes from geometric function theory and nonlinear elasticity. $W$ is polyconvex when restricted to $\text{GL}_+(\R^n)$ but not on the entire $\LL(\R^n)$. Using \eqref{eq:Inv}, we find that 

\begin{align*}
DW(X)&=p\vert X\vert^{p-2}X-p\vert \det(X)\vert\vert X^{-1}\vert^{p-2}(X^\ast X X^\ast)^{-1}+\text{sgn}(\det(X))\vert X^{-1}\vert^{p}\text{adj}(X)^\ast,
\end{align*}

\begin{align*}
X^\ast DW(X)&=p\vert X\vert^{p-2}X^\ast X-p\vert \det(X)\vert\vert X^{-1}\vert^{p-2}X^\ast(X^\ast X X^\ast)^{-1}+\text{sgn}(\det(X))\vert X^{-1}\vert^{p}X^\ast\text{adj}(X)^\ast\\
&=p\vert X\vert^{p-2}X^\ast X-p\vert \det(X)\vert\vert X^{-1}\vert^{p-2}(X^\ast X)^{-1}+\text{sgn}(\det(X))\vert X^{-1}\vert^{p}\det(X)I\\
&=p\vert X\vert^{p-2}X^\ast X-p\vert \det(X)\vert\vert X^{-1}\vert^{p-2}(X^\ast X)^{-1}+\vert \det(X)\vert\vert X^{-1}\vert^{p} I
\end{align*}

Thus
\begin{align*}
T(X)=p\vert X\vert^{p-2}X^\ast X-p\vert \det(X)\vert\vert X^{-1}\vert^{p-2}(X^\ast X)^{-1}-\vert X\vert^{p} I.
\end{align*}

The reduced energy-momentum tensor becomes

\begin{align*}
\mathcal{T}(Y)=p(\text{tr}(Y))^{(p-2)/2}Y-p\sqrt{\det(Y)}\text{tr}(Y^{-1})^{(p-2)/2}Y^{-1}-\text{tr}(Y)^{p/2} I.
\end{align*}

We now assume that $p=n=2$. By the Cayley-Hamilton theorem we have 
\begin{align*}
Y^{-1}=\frac{1}{\det(Y)}[\text{tr}(Y)I-Y].
\end{align*}

This gives 
\begin{align*}
\mathcal{T}(Y)&=2Y-2\sqrt{\det(Y)}Y^{-1}-\text{tr}(Y)I=2Y-2\sqrt{\det(Y)}\frac{1}{\det(Y)}[\text{tr}(Y)I-Y]-\text{tr}(Y)I\\
&=2Y-2\frac{\text{tr}(Y)I}{\sqrt{\det(Y)}}+\frac{2}{\sqrt{\det(Y)}}Y-\text{tr}(Y)I\\
&=2\bigg(1+\frac{1}{\sqrt{\det(Y)}}\bigg)Y-\bigg(1+\frac{2}{\sqrt{\det(Y)}}\bigg)\text{tr}(Y)I
\end{align*}

Let $Y=\alpha I$, $\alpha>0$. Then 
\begin{align*}
\mathcal{T}(\alpha I)&=2\bigg(1+\frac{1}{\alpha}\bigg)\alpha I-\bigg(1+\frac{2}{\alpha}\bigg)\text{tr}(\alpha I)I\\
&=2\bigg(1+\frac{1}{\alpha}\bigg)\alpha I-\bigg(1+\frac{2}{\alpha}\bigg)2\alpha I\\
&=2\alpha\bigg(1+\frac{1}{\alpha}-1-\frac{2}{\alpha}\bigg)I=-2I.
\end{align*}
Thus $\mathcal{T}$ is not injective, but there exists no two $\alpha_1\neq \alpha_2$ so that $\alpha_1I$ and $\alpha_2I$ satisfy the assumptions of Theorem \ref{thm:DW}. 
 
\end{ex}

\begin{ex}[Ball class]
\label{ex:Ball}
Let $\Om \subset \R^n$ be a Lipschitz domain and define the class of mappings
\begin{align*}
\mathscr{A}_{p,q}:=\{u:\Om \to \R^n: Du\in L^p, \,\, \adj(Du)\in L^q\},
\end{align*}
where $p\geq n-1$ and $q\geq p/(p-1)$. This is studied in \cite{Ball3,Ball4, Sv}, see also \cite{FH} for a similar class of polyconvex functionals.
The associated energy to this function class is 
\begin{align*}
I[u]=\int_{\Om}\vert Du(x)\vert^p+\vert \text{adj}(Du)\vert^qdx. 
\end{align*}
and 
\begin{align*}
W(X)=\vert X\vert^p+\vert \text{adj}X\vert^q
\end{align*}
which is a frame indifferent strictly polyconvex function. Using that $W(X)= \vert \tr(X^\ast X)\vert^{p/2}+ \vert \tr(\adj(X)^\ast X)\vert^{q/2}$ and Lemma \ref{lem:MComp} we have 
\begin{align*}
DW(X)&=p\vert X\vert^{p-2}X+q\vert \adj(X)\vert^{q-2}(-\det(X)^{-1}\adj(X^\ast X X^\ast)+\det(X)^{-1}\vert \adj(X)\vert^2\adj(X)^\ast)\\
T(X)&=p\vert X\vert^{p-2}X^\ast X+q\vert \adj(X)\vert^{q-2}(-\adj(X^\ast X)+\vert \adj(X)\vert^2I)
\end{align*}

and thus

\begin{align*}
T(I)&=p\vert I\vert^{p-2}I^\ast I+q\vert \adj(I)\vert^{q-2}(-\adj(I^\ast I)+\vert \adj(I)\vert^2I)\\
&=pn^{(p-2)/2}I+q(n-1)n^{(q-2)/2}I.
\end{align*}
Thus $T$ does not satisfy the assumptions of Proposition \ref{prop:Lam}. 

The reduced energy momentum tensor becomes 
\begin{align*}
\mathcal{T}(Y)&=p\vert \tr(Y)\vert^{(p-2)/2}Y+q\vert \tr(\adj(Y))\vert^{(q-2)/2}(-\adj(Y)+ \tr(\adj(Y)I)
\end{align*}

We now specialise to dimension $n=2$ and choose $p=q=2$. We get 
\begin{align*}
\mathcal{T}(Y)&=2Y+2(-\adj(Y)+\tr(\adj(Y))I).
\end{align*}
If $Y$ is a diagonal matrix then 
$\tr(\adj(Y))=\tr(Y)$ and we get 
\begin{align*}
\mathcal{T}(Y)&=2Y-2\adj(Y)+2\tr(Y)I=2Y
\end{align*}
We now consider the equation $\mathcal{T}(Y)=Z$ for some positive definite diagonal matrix $Z$, which implies that $Y=Z/2$. Thus we have a unique solution and the assumptions of Theorem \ref{thm:DW} are not satisfied. 
\end{ex}

\section*{ Appendix}

\subsection*{ Notions of convexity}

We here recall the notions of polyconvexity and quasiconvexity.

\begin{Def}[Rank-one convexity]
A locally bounded Borel measurable function $W:\LL(\R^n,\R^m)\to \R$ is called \emph{rank-one convex} if for every $A,B\in \LL(\R^n,\R^m)$ such that $\text{rank}(B)\leq 1$ the function $\varphi:\R\to \R\cup\{\infty\}$ given by
\begin{align*}
\varphi(t)=W(A+tB)
\end{align*}
is convex. If in addition for all such $A$ and $B$ the function $\varphi$ is strictly convex we say that $W$ is \emph{strictly rank-one convex}. If $W\in C^2(\LL(\R^n,\R^m))$, then the rank-one convexity is equivalent to the ellipticity condition (also called \emph{Legendre-Hadamard condition})
\begin{align*}
\langle u\otimes v, D^2W(X)(u\otimes v)\rangle \geq 0
\end{align*}
for every $X\in \LL(\R^n,\R^m)$ and every $u\in \R^n$ and $v\in \R^m$. 
\end{Def}

\begin{Def}[Quasiconvex]
A locally bounded Borel measurable function $W:\LL(\R^n,\R^m)\to \R$ is called \emph{quasiconvex} if 
\begin{align*}
W(A)\leq \frac{1}{\vert B_1(0)\vert}\int_{B_1(0)}W(A+D\phi(x))dx
\end{align*}
for every $A\in \LL(\R^n,\R^m)$ and every $\phi\in W^{1,\infty}_0(B_1(0),\R^m)$
\end{Def}

\begin{Def}[Exterior extension of linear map]
\label{def:ExE}
Let $T\in \LL(\R^n,\R^m)$. The \emph{exterior extension} $\Lambda T$ is the unique extension of $T$ to an exterior algebra homomorphism $\Lambda T: \Lambda \R^n \to \Lambda \R^m$ such that 
\begin{itemize}
\item[(i)] $\Lambda T(1)=1$.
\item[(ii)] $\Lambda T(v_1\wedge v_2\wedge...\wedge v_k)=T(v_1)\wedge T(v_2)\wedge ...\wedge T(v_k)$ for any vectors $v_1,v_2,...,v_k\in \R^n$. 
\item[(iii)] $\Lambda T(\alpha w_1+\beta w_2)=\alpha \Lambda T(w_1)+\beta \Lambda T(w_2)$ for any $\alpha,\beta\in \R$ and any $w_1,w_2\in \Lambda \R^n$. 
\end{itemize}
We furthermore let $\Lambda^kT=\Lambda T\vert_{\Lambda^k \R^n}$. 
\end{Def}

Note that in particular $\Lambda^0T(\alpha)=\alpha$ for all $\alpha\in \R$ and $\Lambda^nT(w)=\det(T)w$ for $w\in \Lambda^n\R^n$. Furthermore, if $\text{rank}(T)=k$ then $\Lambda^lT=0$ for $l>k$. 

\begin{Def}[Grade preserving linear map]
A linear map $T\in \LL(\Lambda \R^n,\Lambda\R^m)$ is called \emph{grade preserving} if $T(\Lambda^k\R^n)\subset \Lambda^k\R^m$ for every $k=0,1,2,...,n$. The space of grade preserving linear maps will be denoted by $\widehat{\LL}(\Lambda \R^n,\Lambda \R^m)$. 
\end{Def}

Note that $\text{dim}(\widehat{\LL}(\Lambda \R^n,\Lambda \R^m))=\sum_{k=0}^{\text{min}\{n,m\}}\binom{m}{k}\binom{n}{k}$. In particular, in the case when $m=n$ then $\text{dim}(\widehat{\LL}(\Lambda \R^n,\Lambda \R^n))=\binom{2n}{n}$.

\begin{Def}[Polyconvex]
A locally bounded Borel measurable function $W:\LL(\R^n,\R^m)\to \R$ is called \emph{polyconvex} if there exists a convex function $\sigma: \widehat{\LL}(\Lambda \R^n,\Lambda \R^m)\to \R$ such that 
\begin{align*}
W(X)=\sigma(\Lambda X). 
\end{align*}
We say that $W$ is \emph{strictly polyconvex} if $\sigma$ is strictly polyconvex. 
\end{Def}

In coordinates polyconvexity means that $W$ can be written as a convex function of $X$ and all its minors. Also we typically let $\sigma$ be independent of $\Lambda^0T=\text{id}\vert_{\R}$ so that $W(X)=\sigma(T,\Lambda^2T,\Lambda^3T,...,\Lambda^n T )$.

\subsection*{\sffamily Matrix computations}

For the convenience of the reader we here state a number of useful results regarding functions of matrices.

\begin{Lem}
Let $W(X)=\det(X)$. Then 
\begin{align}\label{eq:Det}
DW(X)=\adj(X)^\ast. 
\end{align}
\end{Lem}

\begin{proof}
For $X,H\in \LL(\R^n)$ we have the expansion
\begin{align*}
\det(X+\eps H)=\det(X)+\eps\text{tr}(\text{adj}(X)H)+O(\eps^2).
\end{align*}

Thus if $W(X)=\det(X)$, then 
\begin{align*}
\langle DW(X),H\rangle &=\lim_{\eps \to 0}\frac{W(X+\eps H)-W(X)}{\eps}=\lim_{\eps \to 0}\frac{\det(X)+\eps\text{tr}(\text{adj}(X)H)+O(\eps^2)-W(X)}{\eps}\\
&=\text{tr}(\text{adj}(X)H)=\langle \text{adj}(X)^\ast,H\rangle.
\end{align*}
\end{proof}

\begin{Lem}
Let $W(X)=\vert X^{-1}\vert^2=\tr((X^{-1})^\ast X^{-1})$. Then for $\det(X)\neq 0$
\begin{align}\label{eq:Inv}
DW(X)=-2(X^{-1})^\ast X^{-1} (X^{-1})^\ast. 
\end{align}
More generally, for $W_p(X)=\vert X^{-1}\vert^p$ and $\det(X)\neq 0$
\begin{align}\label{eq:Inv2}
DW(X)=-p\vert X^{-1}\vert^{p-2}(X^{-1})^\ast X^{-1} (X^{-1})^\ast.
\end{align} 
\end{Lem}

\begin{proof}
For $X,H\in \LL(\R^n)$ we have the expansion
\begin{align*}
(X+\eps H)^{-1}=X^{-1}-\eps X^{-1}HX^{-1}+O(\eps^2).
\end{align*}

Thus, 
\begin{align*}
&((X+\eps H)^{-1})^\ast(X+\eps H)^{-1}=(X^{-1}-\eps X^{-1}HX^{-1}+O(\eps^2))^\ast (X^{-1}-\eps X^{-1}HX^{-1}+O(\eps^2))\\
&=(X^{-1})^\ast X^{-1}-\eps( (X^{-1}HX^{-1})^\ast X^{-1}-(X^{-1})^\ast X^{-1}HX^{-1})+O(\eps^2),
\end{align*}
and 
\begin{align*}
W(X+\eps H)&=\tr((X^{-1})^\ast X^{-1})-\eps \tr((X^{-1}HX^{-1})^\ast X^{-1}-(X^{-1})^\ast X^{-1}HX^{-1})+O(\eps^2)\\
&=W(X)-\eps\tr((X^{-1})^\ast H^\ast(X^{-1})^\ast X^{-1})-\eps \tr((X^{-1})^\ast X^{-1}HX^{-1})+O(\eps^2)\\
&=W(X)-\eps\tr( H^\ast(X^{-1})^\ast X^{-1}(X^{-1})^\ast)-\eps \tr(X^{-1}(X^{-1})^\ast X^{-1}H)+O(\eps^2)\\
&=W(X)-2\eps\tr( H^\ast(X^{-1})^\ast X^{-1}(X^{-1})^\ast)+O(\eps^2)\\
&=W(X)-2\eps \langle (X^{-1})^\ast X^{-1}(X^{-1})^\ast,H\rangle +O(\eps^2).
\end{align*}

\end{proof}

\begin{Lem}
Let $n=2$ and assume that
\begin{align*}
W(X)=\sigma(\tr(X^\ast X),\det(X^\ast X))
\end{align*}
for some smooth function $\sigma:\R^2\to \R$ such that $W$ is smooth. Then, with $Y=X^\ast X$ 
\begin{align*}
T(X)=2\dv_1\sigma(\tr(Y),\det(Y))Y+2\dv_2\sigma(\tr(Y),\det(Y))\det(Y)I-\sigma(\tr(Y),\det(Y))I
\end{align*}
\end{Lem}

\begin{proof}
Since 
\begin{align*}
\tr((X+\eps H)^\ast (X+\eps H))&=\tr(X^\ast X)+2\eps \langle X,H\rangle+O(\eps^2),\\
\det((X+\eps H)^\ast (X+\eps H))&=\det(X^\ast X)+2\eps \langle X\adj(X^\ast X),H\rangle +O(\eps^2)
\end{align*}
we have
\begin{align*}
&W(X+\eps H)=\sigma(\tr(X^\ast X) ,\det(X^\ast X))\\&+2\eps \dv_1\sigma(\tr(X^\ast X) ,\det(X^\ast X))\langle X,H\rangle+2\eps \dv_2\sigma(\tr(X^\ast X) ,\det(X^\ast X))\langle X\adj(X^\ast X),H\rangle+O(\eps^2), 
\end{align*}
and 
\begin{align*}
DW(X)&=2 \dv_1\sigma(\tr(X^\ast X) ,\det(X^\ast X))X+2 \dv_2\sigma(\tr(X^\ast) ,\det(X^\ast X)) X\adj(X^\ast X)\\ 
X^\ast DW(X)&=2 \dv_1\sigma(\tr(X^\ast X) ,\det(X^\ast X))X^\ast X+2 \dv_2\sigma(\tr(X^\ast X) ,\det(X^\ast X)) X^\ast X\adj(X^\ast X)\\
&=2 \dv_1\sigma(\tr(X^\ast X) ,\det(X^\ast X))X^\ast X+2 \dv_2\sigma(\tr(X^\ast X) ,\det(X^\ast X)) X^\ast X\adj(X) \adj(X^\ast)\\
&=2 \dv_1\sigma(\tr(X^\ast X) ,\det(X^\ast X))X^\ast X+2 \dv_2\sigma(\tr(X^\ast X) ,\det(X^\ast X))\det(X^\ast X)I 
\end{align*}

\end{proof}

\begin{Lem}\label{lem:Adj1}
\begin{align*}
\adj(X+\eps H)=\adj(X)-\eps \adj(X)HX^{-1}+\eps\tr(\adj(X) H)X^{-1}+O(\eps^2)
\end{align*}
\end{Lem}

\begin{proof}
If $\det(X)\neq 0$ then 
\begin{align*}
\adj(X+\eps H)&=\det(X+\eps H)(X+\eps H)^{-1}\\&=(\det(X)+\eps \tr(\adj(X) H)+O(\eps^2))(X^{-1}-\eps X^{-1}HX^{-1}+O(\eps^2))\\
&=\det(X)X^{-1}-\eps \det(X)X^{-1}HX^{-1}+\eps \tr(\adj(X) H)X^{-1}+O(\eps^2)\\
&=\adj(X)-\eps \adj(X)HX^{-1}+\eps \tr(\adj(X) H)X^{-1}+O(\eps^2)
\end{align*}

\end{proof}

\begin{Lem}\label{lem:Adj2}
If $W(X)=\tr(\adj(X)^\ast \adj(X))$ and if $\det(X)\neq0$ then
\begin{align*}
DW(X)=-2\det(X)^{-1}\adj(X^\ast X X^\ast)+2\det(X)^{-1}\vert \adj(X)\vert^2\adj(X)^\ast
\end{align*}
and
\begin{align*}
X^\ast DW(X)&=-2\adj(X^\ast X)+2\vert \adj(X)\vert^2I.
\end{align*}
\end{Lem}

\begin{proof}
Using Lemma \ref{lem:Adj1} and the identities $\adj(XY)=\adj(Y)\adj(X)$ and $\adj(X)^\ast=\adj(X^\ast)$ 

\begin{align*}
&\adj(X+\eps H)^\ast \adj(X+\eps H))=(\adj(X)-\eps \adj(X)HX^{-1}+\eps\tr(\adj(X) H)X^{-1}+O(\eps^2))^\ast\\&(\adj(X)-\eps \adj(X)HX^{-1}+\eps\tr(\adj(X) H)X^{-1}+O(\eps^2))\\
&=(\adj(X^\ast)-\eps (X^\ast)^{-1}H^\ast\adj(X^\ast)+\eps\tr(\adj(X) H)(X^\ast)^{-1}+O(\eps^2))\\
&(\adj(X)-\eps \adj(X)HX^{-1}+\eps\tr(\adj(X) H)X^{-1}+O(\eps^2))\\
&=\adj(X^\ast)\adj(X)-\eps \adj(X^\ast)\adj(X)HX^{-1}+\eps\tr(\adj(X) H)\adj(X^\ast)X^{-1}\\
&-\eps (X^\ast)^{-1}H^\ast\adj(X^\ast)\adj(X)+\eps\tr(\adj(X) H)(X^\ast)^{-1}\adj(X)+O(\eps^2))\\
&=\adj(XX^\ast)-\eps \adj(XX^\ast)HX^{-1}+\eps\tr(\adj(X) H)\det(X^{-1})\adj(XX^\ast)\\
&-\eps (X^\ast)^{-1}H^\ast\adj(XX^\ast)+\eps\tr(\adj(X) H)\det(X^{-1})\adj(XX^\ast)+O(\eps^2)).
\end{align*}

Thus,

\begin{align*}
&\tr(\adj(X+\eps H)^\ast \adj(X+\eps H))\\&=\tr(\adj(XX^\ast))-\eps \tr(\adj(XX^\ast)HX^{-1})+\eps\tr(\adj(X) H)\det(X^{-1})\tr(\adj(XX^\ast))\\
&-\eps \tr((X^\ast)^{-1}H^\ast\adj(XX^\ast))+\eps\tr(\adj(X) H)\det(X^{-1})\tr(\adj(XX^\ast))+O(\eps^2))\\
&=\tr(\adj(XX^\ast))-\eps \tr(X^{-1}\adj(XX^\ast)H)+\eps\tr(\adj(X) H)\det(X^{-1})\tr(\adj(XX^\ast))\\
&-\eps \tr(H^\ast\adj(XX^\ast)(X^\ast)^{-1})+\eps\tr(\adj(X) H)\det(X^{-1})\tr(\adj(XX^\ast))+O(\eps^2))\\
&=\tr(\adj(XX^\ast))-\eps \det(X^{-1})\tr(\adj(XX^\ast X)H)+\eps\tr(\adj(X) H)\det(X^{-1})\tr(\adj(XX^\ast))\\
&-\eps \det(X)^{-1}\tr(H^\ast\adj(X^\ast XX^\ast))+\eps\tr(\adj(X) H)\det(X^{-1})\tr(\adj(XX^\ast))+O(\eps^2))\\
&=\tr(\adj(XX^\ast))-\eps \det(X^{-1})\langle (\adj(XX^\ast X))^\ast ,H\rangle +\eps\langle \adj(X)^\ast, H\rangle \det(X^{-1})\tr(\adj(XX^\ast))\\
&-\eps \det(X)^{-1}\langle H,\adj(X^\ast XX^\ast)\rangle+\eps\langle \adj(X)^\ast, H\rangle \det(X^{-1})\tr(\adj(XX^\ast))+O(\eps^2))\\
&=\tr(\adj(XX^\ast))-2\eps \langle \det(X)^{-1}\adj(X^\ast X X^\ast)  ,H\rangle +2\eps\langle \det(X)^{-1}\tr(\adj(XX^\ast))\adj(X)^\ast, H\rangle +O(\eps^2)).
\end{align*}
Hence 
\begin{align*}
DW(X)=-2\det(X)^{-1}\adj(X^\ast X X^\ast)+2\det(X)^{-1}\vert \adj(X)\vert^2\adj(X)^\ast
\end{align*}
and 
\begin{align*}
X^\ast DW(X)&=-2\det(X)^{-1}X^\ast \adj(X^\ast X X^\ast)+2\det(X)^{-1}\vert \adj(X)\vert^2X^\ast \adj(X)^\ast\\&=-2\adj(X^\ast X)+2\vert \adj(X)\vert^2I.
\end{align*}
\end{proof}

Combining the previous Lemmata we find:

\begin{Lem}\label{lem:MComp}
Let $\Phi_i:\{x\in \R: x>0\}\to \R$ be $C^1$ for $i=1,2,3$ and let $W_1(X)=\Phi_1( \tr(X^\ast X))$,  $W_2(X)=\Phi_2( \tr(\adj(X^\ast) \adj(X)))$ and $W_3(X)=\Phi_3( \det(X^\ast X))$ for $X\in \LL(\R^n)$. Then 
\begin{align*}
DW_1(X)&=2\Phi'_1( \tr(X^\ast X))X,\\
DW_2(X)&=2\Phi_2'( \tr(\adj(X^\ast) \adj(X)))(-\det(X)^{-1}\adj(X^\ast X X^\ast)+\det(X)^{-1}\vert \adj(X)\vert^2\adj(X)^\ast),\\
DW_3(X)&=2\Phi_3'( \det(X^\ast X))X\adj(X^\ast X),
\end{align*}
and 
\begin{align*}
X^\ast DW_1(X)&=2\Phi_1'( \tr(X^\ast X))X^\ast X,\\
X^\ast DW_2(X)&=2\Phi_2'(\tr(\adj(X^\ast) \adj(X)))(-\adj(X^\ast X)+\vert \adj(X)\vert^2I),\\
X^\ast DW_3(X)&=2\Phi_3'(  \det(X^\ast X))\det(X^\ast X)I. 
\end{align*}
\end{Lem}

\subsection*{Acknowledgements}
Erik Duse was supported by the Knut and Alice Wallenberg Foundation grant KAW 2015.0270. The author thanks Daniel Faraco for providing references and explaining results regarding differential inclusions into the orthogonal group. Finally, the author thanks Pekka Pankka for interesting discussions on the energy-momentum equations.

{\sc Erik Duse}, Department of Mathematics and Statistics, KTH,  Stockholm, Sweden \texttt{duse@kth.se}

\end{document}